\documentclass{styles/svproc}

\usepackage{url}
\usepackage{amsmath,bm}
\usepackage{amssymb}
\usepackage{tikz}
\usetikzlibrary{matrix,tikzmark,calc,,arrows,shapes,decorations.pathreplacing}
\usepackage{float}
\usepackage{mathtools}
\usepackage{caption} 
\def\curl{\operatorname{curl}}
\def\Curl{\operatorname{Curl}}

\def\dy{\partial_y}

\def\dx{\partial_{x}}

\def\RR{\mathbb{R}}

\def\ttB{\mathsf{B}}
\def\ttC{\mathsf{C}}

\def\ttM{\mathsf{M}}

\def\ttr{\mathsf{r}}

\def\ttw{\mathsf{w}}

\def\tth{\mathsf{h}}
\def\tta{\mathsf{a}}

\def\tth{\mathsf{h}}

\def\ttL{\mathsf{L}}

\def\Th{\mathcal{T}_h}
\def\Eh{\mathcal{E}_h}
\def\N{\mathcal{N}}

\spnewtheorem{mylemma}[lemma]{Lemma}{\bfseries}{\itshape}
\spnewtheorem{myproblem}[lemma]{Problem}{\bfseries}{\itshape}
\spnewtheorem{mytheorem}[lemma]{Theorem}{\bfseries}{\itshape}

\begin{document}
\mainmatter              %
\title{An efficient mixed finite element method for \\ nonlinear magnetostatics and quasistatics}
\titlerunning{Mixed FEM for nonlinear magneto-quasistatics}  %
\author{Herbert Egger\inst{1,2} \and Felix Engertsberger\inst{1} \and Bogdan Radu\inst{2}}
\authorrunning{H. Egger and B. Radu} %
\institute{{Institute~of~Numerical~Mathematics,~{Johannes}~Kepler~University,~Linz,~Austria} \\
\and
{Johann~Radon~Institute~for~Computational~and~Applied~Mathematics,~Linz,~Austria} \\
\email{herbert.egger@jku.at}, 
\email{felix.engertsberger@jku.at}, \email{bogdan.radu@ricam.oeaw.ac.at}}

\maketitle %

\begin{abstract}
We consider systems of nonlinear magnetostatics and quasistatics that typically arise in the modeling and simulation of electric machines. The nonlinear problems, eventually obtained after time discretization, are usually solved by employing a vector potential formulation. In the relevant two-dimensional setting, a discretization can be obtained by $H^1$-conforming finite elements. 
We here consider an alternative formulation based on the $H$-field which leads to a nonlinear saddlepoint problem. After commenting on the unique solvability, we study the numerical approximation by $H(\curl)$-conforming finite elements and present the main convergence results. 
A particular focus is put on the efficient solution of the linearized systems arising in every step of the nonlinear Newton solver. 
Via hybridization, the linearized saddlepoint systems can be transformed into linear elliptic problems, which can be solved with similar computational complexity as those arising in the vector or scalar potential formulation. 
In summary, we can thus claim that the mixed finite element approach based on the $H$-field can be considered a competitive alternative to the standard vector or scalar potential formulations for the solution of problems in nonlinear magneto-quasistatics. 

\end{abstract}

\section{Introduction}
We consider nonlinear systems of magnetostatics and quasistatics which typically arise in the modeling of electric machines. After eventual time discretization, the governing equations in the relevant two-dimensional setting can be phrased as \begin{alignat}{2}
g'(H) - \Curl a &= 0 \qquad &&  \text{in } \Omega, \label{eq:1}\\
\curl H + \sigma a &= j \qquad && \text{in } \Omega, \label{eq:2}
\end{alignat}
and we impose homogeneous boundary conditions $a=0$ on $\partial\Omega$ 
for simplicity. 
Here $\Omega \subset \RR^2$ is the cross-section of the machine, $H$ the parallel component of the magnetic field intensity, $a$ the axial component of the magnetic vector potential, and $\curl H= \dy H_x - \dx H_y$ and $\Curl a=(\dy a, -\dx a)$ are the two curl operators arising in two space dimensions. 
Furthermore, $g'(H) = \nabla_H g(H)$ is the derivative of the magnetic coenergy density $g$ which, in general, may depend on the spatial coordinate as well. 
Finally, $j$, $\sigma$ are the axial component of the current density and a  generalized conductivity, with $\sigma=0$ in the static case.
Before we proceed, let us briefly discuss some of the standard approaches for the analysis and numerical solution of the problem~\cite{Henrotte97,meunier08}.

\subsection*{Magnetic vector potential formulation}
We can use \eqref{eq:1} to express $H = f'(\Curl a)$, where $f'(\cdot)$ is the inverse of $g'(\cdot)$, i.e., the derivative of the magnetic energy density~\cite{silvester91}. In fact, $f(\cdot)$ and $g(\cdot)$ are closely linked via convex duality~\cite{rockafellar70}. 
Inserting this into \eqref{eq:2} leads to 
\begin{align} \label{eq:3}
\sigma a + \curl (f'(\Curl a)) &= j \qquad \text{in } \Omega.
\end{align}
Any solution to this problem with zero boundary conditions can be characterized by the following weak formulation: Find $u \in H_0(\Curl;\Omega)$ such that 
\begin{align} \label{eq:4}
(\sigma a, v) + (f'(\Curl a), \Curl v) &= (j,v) \qquad \forall v \in H_0(\Curl;\Omega).
\end{align}
Since $\Curl a = (\nabla a)^\perp$, one can see that $H_0(\Curl;\Omega) \simeq H_0^1(\Omega)$ are isomorphic. 
Under some general assumptions, the existence of a unique solution can then be proven by a nonlinear version of the Lax-Milgram theorem~\cite{Zeidler1990}, and standard finite elements can be used to approximate \eqref{eq:4} numerically; see e.g.~\cite{heise94,meunier08}. 

\subsection*{Total magnetic field formulation}

If $\sigma>0$ everywhere, one can use \eqref{eq:2} to express $a = \rho (j - \curl H)$ with modified resistivity $\rho=\frac{1}{\sigma}$. This leads to the pure $H$-field formulation
\begin{align} \label{q:5}
 g'(H) + \Curl(\rho \curl H) &= \Curl(\rho \,  j) \qquad \text{in } \Omega.
\end{align}
The weak form of this problem reads: Find $H \in H(\curl;\Omega)$ such that 
\begin{align} \label{eq:6}
(g'(H), w) + (\rho \curl H, \curl w) &= (\rho \, j, \curl w) \qquad \forall w \in H(\curl;\Omega). 
\end{align}
Note that the space $H(\curl;\Omega)$ consists of vector fields with tangential continuity. The existence of a unique solution to this nonlinear elliptic problem can be established with similar arguments as before.
Furthermore, N\'ed\'elec finite elements can be used for the systematic numerical approximation~\cite{Bossavit88,Nedelec80}. 

\subsection*{Reduced scalar potential formulation}

In the static case, we have $\sigma = 0$ and one can satisfy \eqref{eq:2} explicitly by decomposing $H = H_s + \nabla \psi$, with given source field $H_s$ such that $\curl H_s = j$. By testing equation~\eqref{eq:2} with $\nabla z$, one arrives at the following weak formulation: Find $\psi \in H^1(\Omega) / \RR$ such that 
\begin{align}
(g'(H_s + \nabla \psi), \nabla z) &= 0 \qquad \forall z \in H^1(\Omega)/\RR. 
\end{align}
Here $H^1(\Omega)/\RR$ denotes the space of scalar potentials with zero average. 
This is again a nonlinear elliptic problem and the existence of a unique solution can be established as before~\cite{Engertsberger23,Zeidler1990}. Moreover, standard $H^1$-conforming finite elements can be used for the efficient and reliable numerical approximation \cite{Albanese90,Engertsberger23}.

\subsection*{Outline and main contributions}

In this paper, we consider the mixed formulation \eqref{eq:1}--\eqref{eq:2} of magneto-quasistatics with $H$ in $H(\curl;\Omega)$ and $a \in L^2(\Omega)$, which works for all $\sigma \ge 0$.
The weak form of this system leads to a nonlinear saddlepoint problem, and we briefly discuss its well-posedness.
Mixed finite elements can be used for the systematic discretization; see e.g. \cite{BarbaMariniSavini1993,RogierSegre92,GilletteBajaj2011}, where corresponding linear magnetostatic problems were considered.
For the solution of the nonlinear problem, we apply a damped Newton method, and we discuss the efficient solution of the linearized saddlepoint problems arising in every Newton through hybridization. By this one can reduce the linearized problems to symmetric positive definite systems which can be solved efficiently with a computational complexity similar to that of the alternative approaches discussed above. 
A detailed comparison with the vector potential formulation, concerning accuracy and computational efficiency, will be given below for a typical test problem arising in electric machine simulation.

\section{Dual mixed formulation}
Let $\Omega\subset \RR^2$ be a simply connected bounded Lipschitz domain. To avoid technicalities, we further assume that $\partial\Omega$ is polygonal.
The weak formulation of the system \eqref{eq:1}--\eqref{eq:2} with $a=0$ on $\partial\Omega$ can then be stated as follows.
\begin{myproblem}\label{prob:main}
Find $(H,a)\in H(\curl,\Omega)\times L^2(\Omega)$ such that 
\begin{alignat*}{2}
(g'(H),v)-(a,\curl v) &= 0     \qquad &&\forall v\in H(\curl;\Omega) \\
(\curl H,q) + (\sigma a, q)          &= (j,q) \qquad &&\forall q\in L^2(\Omega).
\end{alignat*}
\end{myproblem}
With the usual arguments, one can see that any smooth solution of \eqref{eq:1}--\eqref{eq:2} with homogeneous boundary conditions also solves Problem~\ref{prob:main}.
The following result clarifies the well-posedness of the nonlinear variational problem above.
\begin{mylemma} \label{lem:mixed}
Let $g'\in C^1(\RR^2)$ be uniformly monotone and Lipschitz, i.e., \\[1ex]
(i) \quad $\langle g'(y)-g'(z),y-z \rangle\ge \alpha |y-z|^2$ for all $y,z \in \RR^2$,\\[1ex]
(ii) \quad  $|g'(y)-g'(z)| \le C_a |y-z|$ for all $y,z \in \RR^2$, \\[1ex]
with $C_a,\alpha>0$. 
Then for any $j \in L^2(\Omega)$, the variational Problem~\ref{prob:main} admits a unique solution $(H,a) \in H(\curl,\Omega) \times L^2(\Omega)$, which can be bounded by
\begin{align*}
\|H\|_{H(\curl;\Omega)} + \|a\|_{L^2(\Omega)} \le C\| j\|_{L^2(\Omega)}.
\end{align*}
The constant $C$ only depends on $\alpha$ and $C_a$, and on the domain $\Omega$.
\end{mylemma}
The proof of this result can be accomplished using the surjectivity of the vector-to-scalar $\curl$ operator, an extension of the Brezzi lemma~\cite{Boffi13}, and
the Zarantonello theorem for monotone operators~\cite{Zeidler1990}; details are left to the reader.
Let us note that $g$ may additionally depend on the spatial coordinate, in which case the above properties are assumed to hold uniformly for all $x \in \Omega$. 

\section{Finite element discretization}
Let $\Th$ be a geometrically conforming and shape-regular triangulation of the domain $\Omega$.
By $h=\max_T h_T$ we denote the global mesh size, i.e., the maximal diameter of an element in $\Th$. 
We define the finite element  spaces
\begin{align*}
V_h  = \N_k(\Th) \cap H(\curl,\Omega)  \qquad \text{and} \qquad 
Q_h = P_k(\Th),
\end{align*}
where $X(\Th) = \{v  : v|_T \in X(T)\}$ is used to denote functions which are members of $X(T)$ on every element $T$. Further $\N_k(T) = P_k(T)^2 \oplus \vec x\times \mathring P_k(T)$ is the Nedelec finite elements, $P_k(T)$ are the polynomials of degree $\le k$, and $\mathring P_k(T)$ contains the corresponding homogeneous polynomials; see \cite{Boffi13,Nedelec80} for details.
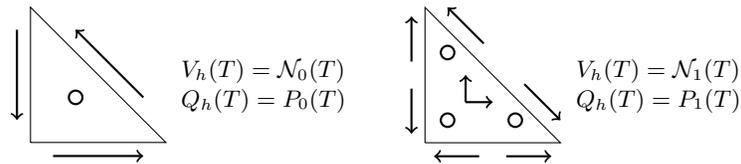
\begin{figure}
    \centering
    \begin{tabular}{llll}
    \begin{tabular}{l}\begin{tikzpicture}[scale=0.6]
\draw (0,0) -- (0,3) -- (3,0) -- (0,0);
\draw[thick,->] (0.5,-0.3) -- (2.5,-0.3);
\draw[thick,->] (-0.3,2.5) -- (-0.3,0.5);
\draw[thick,->] (2.5,1.0) -- (1.0,2.5);
\draw[thick] (1,1) circle (0.15cm);
\end{tikzpicture}\end{tabular}
    & \begin{tabular}{l} $ V_h(T)=\N_0(T)$ \\ $Q_h(T)=P_0(T)$ \end{tabular} \qquad %
    \begin{tabular}{l}\begin{tikzpicture}[scale=0.6]
\draw (0,0) -- (0,3) -- (3,0) -- (0,0);
\draw[thick,->] (1.2,-0.3) -- (0.2,-0.3);
\draw[thick,->] (1.8,-0.3) -- (2.8,-0.3);
\draw[thick,->] (-0.3,1.2) -- (-0.3,0.2);
\draw[thick,->] (-0.3,1.8) -- (-0.3,2.8);
\draw[thick,->] (2.2,1.3) -- (3.0,0.5);
\draw[thick,->] (1.3,2.2) -- (0.5,3.0);
\draw[thick,->] (0.9,0.9) -- (1.5,0.9);
\draw[thick,->] (0.9,0.9) -- (0.9,1.5);

\draw[thick] (2,0.5) circle (0.15cm);
\draw[thick] (0.5,2) circle (0.15cm);
\draw[thick] (0.5,0.5) circle (0.15cm);
\end{tikzpicture}\end{tabular}
    & \begin{tabular}{l} $ V_h(T)=\N_1(T)$ \\ $Q_h(T)=P_1(T)$ \end{tabular} \\ %
    \end{tabular}
    \caption{Degrees of freedom for $V_h$ (arrows) and $Q_h$ (circles) for order $k=0,1$.\label{fig:dofs}}
\end{figure}
For illustration, the degrees of freedom for the spaces of order $k=0$ and $k=1$ are depicted in Figure~\ref{fig:dofs}.
To facilitate the implementation of the nonlinear term, we further introduce an approximation  $(\cdot,\cdot)_h$ for the $L^2$-scalar product $(\cdot,\cdot)$, which is obtained by appropriate numerical integration on each element.
We then consider the following discretization of Problem~\ref{prob:main} by an inexact Galerkin method.
\begin{myproblem}\label{prob:disc}
Find $(H_h,a_h)\in V_h\times Q_h$ such that 
\begin{alignat*}{3}
(g'(H_h),v_h)_h-(a_h,\curl v_h) &= 0     \qquad &&\forall v_h\in V_h \\
(\curl H_h,q_h) + (\sigma a_h,q_h)          &= (j,q_h) \qquad &&\forall q_h\in Q_h.
\end{alignat*}
\end{myproblem}
With similar reasoning as for the continuous problem and standard arguments for the analysis of finite element methods, one can show the following result. 

\begin{mytheorem}
Let $\Omega$, $g'(\cdot)$ and $j$ satisfy the assumption of Lemma~\ref{lem:mixed}, and let $V_h$, $Q_h$ be chosen as above.
Further, let the element quadrature rule underlying $(\cdot,\cdot)_h$ be exact at least for polynomials of degree $ \le 2k+2$ and have strictly positive weights.
Then Problem~\ref{prob:disc} has a unique solution $(H_h,a_h) \in V_h \times Q_h$. 
Moreover
\begin{align*}
\|H-H_h\|_{H(\curl;\Omega)} + \|a-a_h\|_{L^2(\Omega)} \le C h^k (\|H\|_{H^k(\curl;\Th)} + \|a\|_{H^k(\Th)})
\end{align*}
with a constant $C \ne C(h)$ independent of the meshsize $h$. 
\end{mytheorem}
The space $H^k(\curl;\Th)$ consists of functions $v \in H^k(\Th)^2$ with $\curl v \in H_k(\Th)$. 
We refer to \cite{Engertsberger23} and \cite{heise94} for detailed proofs of corresponding results for the vector potential and the scalar potential formulations. With minor modifications, the main arguments can be generalized to the problem under consideration. 

\section{Efficient algebraic solution}
The numerical solution of Problem~\ref{prob:disc} requires a nonlinear solver. We here consider iterative methods of the form a damped Newton iteration, which takes the form 
\begin{align}\label{eq:update}
(H_h^{n+1},a_h^{n+1}) = (H_h^{n},a_h^{n}) + \tau^n (\delta H_h^n,\delta a_h^n), \qquad n \ge 0.
\end{align}
The initial iterate $(H_h^0,a_h^0) \in V_h \times Q_h$ is given. 
The updates $(\delta H_h^n,\delta a_h^n)$ for the Newton iteration are defined by the following \textit{linearized} system.
\begin{myproblem}\label{prob:iter}
Find $(\delta H_h,\delta a_h)\in V_h\times Q_h$ such that 
\begin{alignat}{3}
(g''(H_h^n)\, \delta H_h^{n},v_h)_h-(\delta a_h^{n},\curl v_h) &= -(g'(H_h^n),v_h)_k + (\delta a_h^{n},\curl v_h) \label{eq:lin1}\\
(\curl \delta H_h^{n},q_h) + (\sigma \delta a_h^n, q_h) &= (j,q_h) - (\curl  H_h^{n},q_h) - (\sigma a_h^n,q_h) \label{eq:lin2}
\end{alignat}
for all test functions $v_h\in V_h$ and $q_h\in Q_h$.
\end{myproblem}
Under the previous assumptions, the existence of a unique solution can be obtained here by standard saddlepoint theory~\cite{Boffi13}. 
With an appropriate choice of the step size $\tau^n$, the convergence of the damped Newton iteration \eqref{eq:update} can be established; we refer to \cite{Engertsberger23} for a related analysis.

\subsection*{Hybridization}

The realization of the damped Newton iteration \eqref{eq:update} requires the solution of a linear saddlepoint problem \eqref{eq:lin1}--\eqref{eq:lin2} in every step. This can be achieved efficiently using hybridization \cite{Boffi13,meunier08}.
The key idea is to impose the tangential continuity of the functions in $V_h$ via Lagrange multipliers.
Let $\Eh=\{E_{ij} : i < j\}$ denote the set of all edges $E_{ij} = \partial T_i \cap \partial T_j$ in $\Th$, and consider the spaces
\begin{align*}
\widetilde V_{h} = \N_k(\Th),\qquad Q_h = P_k(\Th),\qquad \widehat Q_h = P_k(\Eh)
\end{align*}
For a piecewise smooth function on $\Th$, set  $[\ell]_{E} := \ell|_{E\cap T_i}+\ell|_{E\cap T_j}$ for any face $E \in \Th$. 
Then Problem~\ref{prob:iter} is equivalent to the following system.
\begin{myproblem}\label{prob:hybrid}
Find $(\delta H_h^{n},\delta a_h^{n},\delta \widehat a_h^{n})\in \widetilde V_{h}\times Q_h\times \widehat Q_h$ such that
\begin{alignat*}{3}
(g''(H_h^n) \, \delta H_h^{n},v_h)_{h,T}&-(\delta a_h^{n},\curl v_h)_T - (\delta \widehat a_h^{n},[n\times v_h])_E  \\
&= -(g'(H_h^n),v_h)_{k,T}+(\delta a_h^{n},\curl v_h)_T + (\delta \widehat a_h^{n},[n\times v_h])_E \\
(\curl \delta H_h^{n} + \sigma \delta a_h^n,q_h)_T             &= (j - \curl H_h^{n} - \sigma a_h^n,q_h)_T \\
([n\times \delta H_h^{n}],\hat q_h)_E      &=0
\end{alignat*}
for all $T\in\Th$ and for all $v_h\in \widetilde V_h$, $q_h\in Q_h$ and $\hat q_h\in \widehat Q_h$. %
\end{myproblem}
Let us note that the ters $[n \times v_h]$ represent the jump of the tangential component of the piecewise smooth vector-valued function $v_h \in \widetilde V_h$. 
The last equation thus ensures continuity of the tangential component of $\delta H_h^n$.
Moreover, one can show that the first two components of the solution of Problem~\ref{prob:hybrid} coincide with the unique solution of Problem~\ref{prob:iter}; we refer to \cite{Boffi13} for details.

\subsection*{Algebraic reduction}

With the usual choice of a basis, the implementation of the mixed finite element approximation of Problem~\ref{prob:hybrid} leads to an algebraic system of the form
\begin{align*}
\begin{pmatrix}
\ttM & -\ttB^\top & -\ttL^\top\\
\ttB & \ttC & 0 \\
\ttL & 0 & 0
\end{pmatrix}
\begin{pmatrix}
\delta\tth\\
\delta\tta\\
\delta\widehat\tta
\end{pmatrix}
=
\begin{pmatrix}
\ttr\\
\ttw\\
0
\end{pmatrix}.
\end{align*}
Since the only coupling between elements occurs via the jump terms in Problem~\ref{prob:hybrid}, one can see that  $\ttM$, $\ttC$, and $\ttB$ are block diagonal, with each block associated with a single element. 
This allows the elimination of both $\delta \tth$ and $\delta \tta$ locally, leading to a Schur complement system for the hybrid variable $\delta \widehat \tta$ alone, viz.
\begin{align}\label{eq:sysmix}
\begin{pmatrix} \ttL \quad  0 \end{pmatrix} 
\begin{pmatrix} \ttM & -\ttB^\top \\ -\ttB & -\ttC \end{pmatrix}^{-1}
\begin{pmatrix} \ttL^\top \\ 0 \end{pmatrix} 
\delta \widehat \tta
&=
\begin{pmatrix} \ttL \quad  0 \end{pmatrix} 
\begin{pmatrix} \ttM & -\ttB^\top \\ -\ttB & -\ttC \end{pmatrix}^{-1}
\begin{pmatrix} -\ttr \\ \ttw \end{pmatrix}
\end{align}
From the previous considerations, one can infer that the system matrix of this problem is sparse, symmetric, and positive definite, and hence \eqref{eq:sysmix} can be solved efficiently. 
Once the solution component $\delta \hat \tta$ at the edges has been computed, the other components $\delta \tth, \delta \tta$ can be obtained by local post-processing; see \cite{Boffi13} for details. 
Let us emphasize that all the required computations can be done efficiently on the algebraic level.
Via hybridization, the computation of the Newton step in Problem~\ref{prob:iter} can thus be achieved by the solution of a global symmetric positive definite system and some local algebraic computations.

\section{Numerical illustration}

To demonstrate the efficiency and accuracy of the proposed hybridized mixed finite element method for nonlinear magnetostatics, let us briefly report on some numerical results for a typical scenario arising in electric machine computation. 
For ease of presentation, we only consider the magnetostatic case $\sigma=0$ here, which is automatically covered by our results.

\subsection*{Model problem}
The geometry, which is taken from~\cite{Gangl2023}, comprises $16$ permanent magnets (PMs) and $48$ slots. 
A linear material law $B=\mu_0 H$, corresponding to $f(B)=\frac{1}{2\mu_0}|B|^2$ resp. $g(H)=\frac{\mu_0}{2} |H|^2$ is used for the slots, the air gap between the stator and rotor, the end pockets of the PMs, and the shaft. 
For the permanent magnets, we use $B=\mu_0(H+M)$ which amounts to $f(B)=\frac{1}{2\mu_0} |B - \mu_0 M|^2$, respectively, and $g(H)=\frac{\mu_0}{2} |H+M|^2$.
The core material of the stator and rotor is described by an isotropic model $f(B)=\tilde f(|B|)$.
For our numerical tests, $\tilde f(|B|)$ is chosen as a cubic spline approximation of the Brauer model~\cite{Brauer75}, which allows us to guarantee all the required properties; see~\cite{Pechstein06} for details.
By convex duality, the coenergy functional then takes the form $g(H)=\tilde g(|H|)$ with $\frac{\tilde g'(|H|)}{|H|} = \frac{|B|}{\tilde f'(|B|)}$. The function $\tilde g(|H|)$ is again approximated by a cubic spline. 
Further details of the geometry and material model can be found in~\cite{Bogdan2023}.
We consider the fields generated solely by the permanent magnets, and therefore choose $j=0$ for the current density. 

\subsection*{Details on the discretization}
In our numerical tests, we consider the first and second-order approximations of the vector potential formulation and the mixed $H$-field model. 
For the implementation of the nonlinear terms, we use numerical quadrature with the Dunavant rules~\cite{Dunavant85} of order $p=2$ and $p=4$, respectively.
Both these rules satisfy our condition on the positivity of the weights used in our theorems.
For both formulations, Armijo backtracking~\cite{Nocedal} is chosen for determining the step size $\tau^n$, which allows to prove global convergence of the damped Newton method for the problems under consideration.

\subsection*{Numerical results}

\begin{figure}
\centering
\includegraphics[scale=0.18]{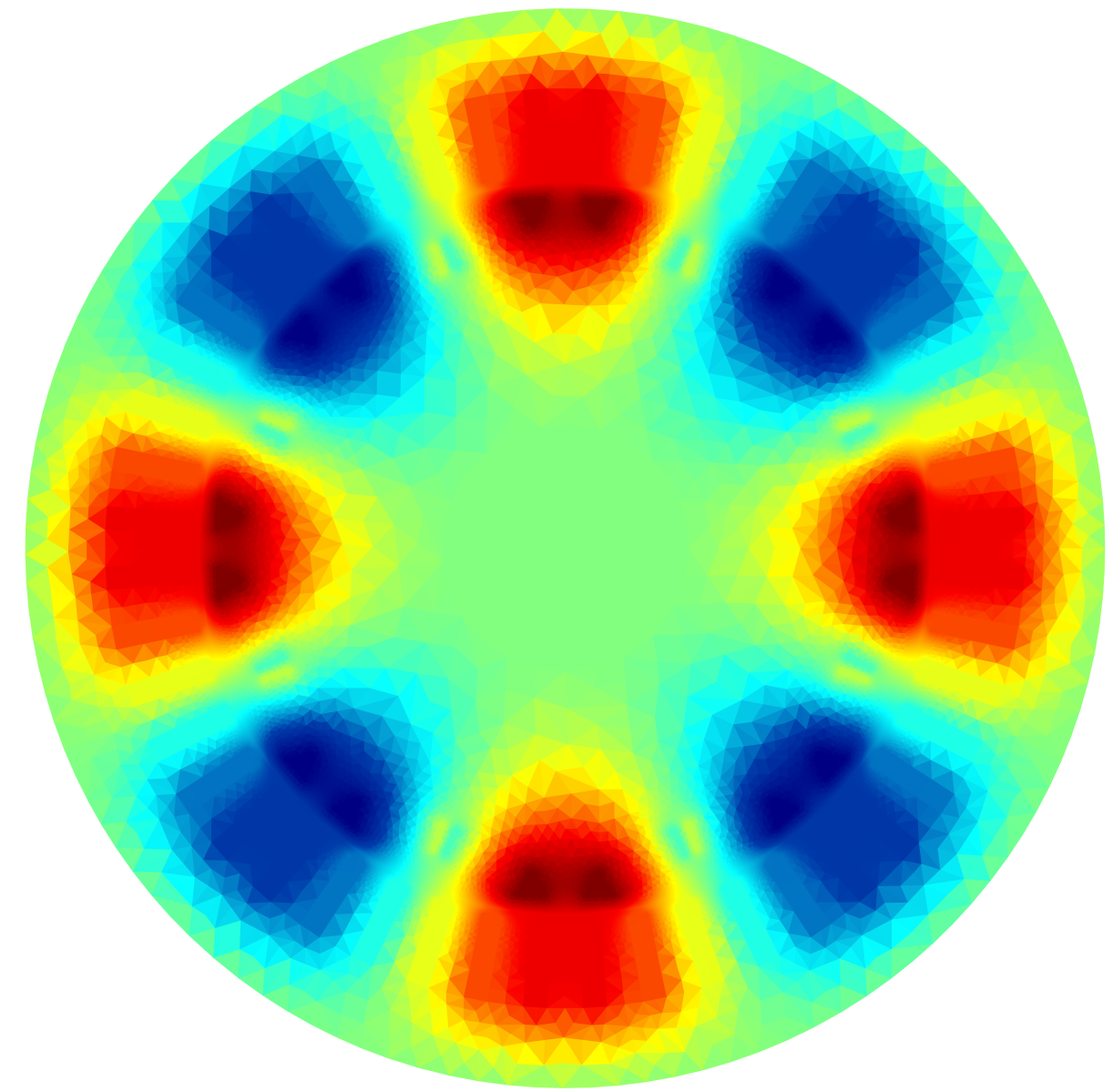}\qquad
\includegraphics[scale=0.18]{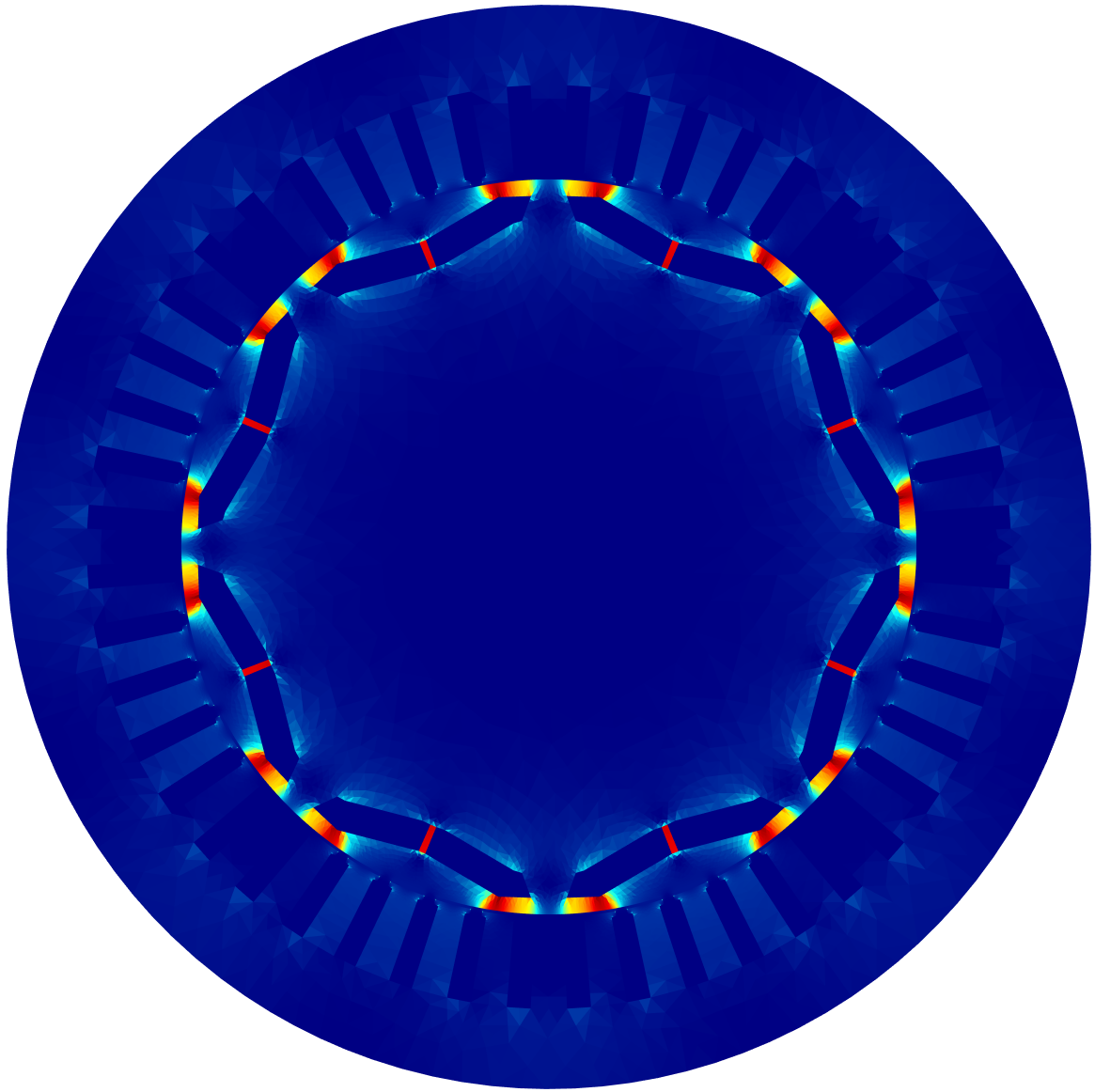}
\caption{Numerical solution of the magnetostatic problem obtained by the mixed $H$-field formulation using lowest finite elements. Left: vector potential $a_h \in P_0(\Th)$; right: magnitude of the magnetic field $H_h \in \N_0(\Th)$.}
\label{fig:solution}
\end{figure}

A typical result of our simulation is depicted in Figure~\ref{fig:solution}.
In the following discussion, we compare the performance and accuracy of the primal and hybridized mixed finite element methods.

Let us next compare the accuracy and computational efficiency of the vector potential and the mixed $H$-field formulation.
Since the exact solution exhibits singularities at the material jumps, only a reduced convergence order can be expected. 
In Table~\ref{tab:errors}, we therefore only report about convergence rates for the lowest order methods. 
\begin{table}[ht!]
\centering
\setlength\tabcolsep{1.5ex}
\renewcommand{\arraystretch}{1.1}
\begin{tabular}{c||c|c|c|c||c|c|c|c} 
$h$ & error & eoc & iter & time 
    & error & eoc & iter & time \\
\hline
\hline
$2^{-1}$ & $0.2938907$ & $-$    & $9$  & $0.005s$
         & $0.306129$ & $-$    & $9$  & $0.01s$ \\
$2^{-2}$ & $0.1745872$ & $0.75$ & $12$ & $0.03s$
         & $0.185483$ & $0.72$ & $11$ & $0.06s$ \\
$2^{-3}$ & $0.1005515$ & $0.79$ & $11$ & $0.16s$
         & $0.110442$ & $0.75$ & $10$ & $0.27s$ \\
$2^{-4}$ & $0.0571080$ & $0.82$ & $13$ & $1.00s$
         & $0.064842$ & $0.77$ & $14$ & $1.60s$ \\
$2^{-5}$ & $0.0323120$ & $0.82$ & $16$ & $8.00s$
         & $0.037660$ & $0.78$ & $15$ & $11.80s$ \\ 
$2^{-6}$ & $-$ & $-$ & $16$ & $66.85$
         & $-$ & $-$ & $18$ & $88.12s$
\end{tabular}
\caption{Results for the lowest order versions of the vector potential formulation (left: $B_h = \Curl a_h$) and the mixed $H$-field formulation (right: $B_h = \pi_h f'(H_h)$). \label{tab:errors} }
\end{table}
In all computations, the damped Newton method converged in about $10-20$ iterations to a relative residual of $10^{-8}$. Both approaches yield similar errors and convergence rates and require similar computation times.

In Table~\ref{fig:comparison}, we report in more detail about the size of the linear systems, the number of non-zero entries, and the solution time required for computing one Newton step. 
We do this for both formulations and the first and second-order methods. 
\begin{table}
\centering
\setlength\tabcolsep{1.5ex}
\renewcommand{\arraystretch}{1.1}
\begin{tabular}{c|c|c|c|c}
method & order & ndofs & nnz & cpu time\\
\hline
\hline
primal & $1$ & $422.015$    & $2.953.627$ & $0.45s$ \\
dual   & $1$ & $1.265.570$  & $6.327.378$ & $1.00s$ \\
\hline
primal & $2$ & $1.687.821$  & $19.368.652$ & $2.30s$ \\
dual   & $2$ & $2.531.140$  & $25.309.512$ & $3.20s$
\end{tabular}
\caption{System size (ndofs), number of non-zeros (nnz), and computation times for sparse solution of one Newton-step obtained in the primal finite element method and the hybrid dual mixed method after elimination of the local dofs. All computations were performed using the \texttt{CHOLMOD} package \cite{CHOLMOD} form SparseSuite on a mesh with $n_p = 422.015$ grid points.}
\label{fig:comparison}
\end{table}
Let us mention that we obtained similar convergence rates for all choices of methods and polynomial degrees, which is to be expected from the theoretical results. 
As can be deduced from the table, hybridization allows to compute the approximation of the mixed $H$-field formulation with similar computational complexity as that of the vector potential approach. For higher polynomial order, the difference in computation times even becomes smaller.

\subsection*{Summary}

In this paper, we considered a mixed finite element approach for solving the magneto-quasistatic problem.
The approach is viable for the magnetostatic and the magneto-quasistatic case, and a damped Newton method can be used to solve the nonlinear problem.
Using hybridization, the Newton-systems can be reduced to scalar elliptic problems, whose numerical solution has a similar efficiency as that of the magnetic scalar or vector potential formulation.
In contrast to the scalar potential approach, the mixed formulation is applicable also in the quasistatic case, and in contrast to the vector potential formulation, it offers direct access to the magnetic field intensity. 
Further information about the magnetic flux could be obtained by post-processing, which also offers access to a-posteriori error estimates; we refer to \cite{Boffi13} for details.

\subsection*{Acknowledgements}
This work was supported by the joint DFG/FWF Collaborative Research Centre CREATOR: Computational Electric Machine Laboratory (TRR361/SFB-F90).


\begin{thebibliography}{10}
\providecommand{\url}[1]{{#1}}
\providecommand{\urlprefix}{URL }
\expandafter\ifx\csname urlstyle\endcsname\relax
  \providecommand{\doi}[1]{DOI~\discretionary{}{}{}#1}\else
  \providecommand{\doi}{DOI~\discretionary{}{}{}\begingroup
  \urlstyle{rm}\Url}\fi

\bibitem{Albanese90}
Albanese, R., Rubinacci, G.: Magnetostatic field computation in terms of
  two-component vector potentials.
\newblock Int. J. Numer. Meth. Engrg. \textbf{29}, 515--532 (1990)

\bibitem{Boffi13}
Boffi, D., Brezzi, F., Fortin, M.: Mixed finite element methods and
  applications.
\newblock Springer, Heidelberg (2013)

\bibitem{Bossavit88}
Bossavit, A.: A rationale for edge-elements in {3-D} fields computations.
\newblock IEEE Trans. Magn. \textbf{24}, 74--79 (1988)

\bibitem{Brauer75}
Brauer, J.: Simple equations for the magnetization and reluctivity curves of
  steel.
\newblock IEEE Trans. Magn. \textbf{11}, 81--81 (1975)

\bibitem{CHOLMOD}
Chen, Y., Davis, T.A., Hager, W.W., Rajamanickam, S.: Algorithm 887: Cholmod,
  supernodal sparse cholesky factorization and update/downdate.
\newblock ACM Trans. Math. Softw. \textbf{35}, art.no. 22 (2008)

\bibitem{BarbaMariniSavini1993}
Di~Barba, P., Marini, L.D., Savini, A.: Mixed finite elements in
  magnetostatics.
\newblock COMPEL \textbf{12}, 113--124 (1993)

\bibitem{Henrotte97}
Dular, P., Remacle, J.F., Henrotte, F., Genon, A., Legros, W.: Magnetostatic
  and magnetodynamic mixed formulations compared with conventional
  formulations.
\newblock IEEE Trans. Magn. \textbf{33}, 1302--1305 (1997)

\bibitem{Dunavant85}
Dunavant, D.A.: High degree efficient symmetrical {G}aussian quadrature rules
  for the triangle.
\newblock Int. J. Numer. Meth. Engrg. \textbf{21}, 1129--1148 (1985)

\bibitem{Engertsberger23}
Engertsberger, F.: The scalar potential approach in nonlinear magnetostatics.
\newblock Master Thesis, Johannes Kepler University Linz, 2023

\bibitem{Gangl2023}
Gangl, P., Gobrial, M., Steinbach, O.: A space-time finite element method for
  the eddy current approximation of rotating electric machines.
\newblock arXiv:2307.00278

\bibitem{GilletteBajaj2011}
Gillette, A., Bajaj, C.: Dual formulations of mixed finite element methods with
  applications.
\newblock Computer-Aided Design \textbf{43}(10), 1213--1221 (2011).
\newblock \doi{https://doi.org/10.1016/j.cad.2011.06.017}.
\newblock
  \urlprefix\url{https://www.sciencedirect.com/science/article/pii/S001044851100159X}.
\newblock Solid and Physical Modeling 2010

\bibitem{heise94}
Heise, B.: Analysis of a fully discrete finite element method for a nonlinear
  magnetic field problem.
\newblock SIAM J. Numer. Anal. \textbf{31}, 745--759 (1994)

\bibitem{meunier08}
Meunier, G.: The Finite Element Method for Electromagnetic Modeling.
\newblock Wiley-ISTE (2008)

\bibitem{Nedelec80}
Nedelec, J.C.: Mixed finite elements in $\mathbb{R}^3$.
\newblock Numer. Math. \textbf{35}, 315--341 (1980)

\bibitem{Nocedal}
Nocedal, J., Wright, S.: Numerical Optimization. 2nd ed.
\newblock Springer (2006)

\bibitem{Pechstein06}
Pechstein, C., Jüttler, B.: Monotonicity-preserving interproximation for {BH}
  curves.
\newblock J. Comput. Appl. Math. \textbf{196}, 45--57 (2006)

\bibitem{Bogdan2023}
Radu, B.: Seminar on electrical machine simulation, {JKU Linz, 2023}.
\newblock $ $\\
  {\small\url{github.com/radu-bogdan/SeminarEM/wiki/Solving-nonlinear-magnetostatics}}

\bibitem{rockafellar70}
Rockafellar, R.T.: Convex Analysis.
\newblock Princeton Univ. Press, Princeton (1970)

\bibitem{RogierSegre92}
Rogier, F., Segre, J.: Mixed finite element method applied to a magnetostatic
  problem.
\newblock Comput. Methods Appl. Mech. Engrg. \textbf{94}, 1--11 (1992)

\bibitem{silvester91}
Silvester, P.P., Gupta, R.P.: Effective computational models for anisotropic
  soft {B-H} curves.
\newblock IEEE Trans. Magn. \textbf{27}, 3804--07 (1991)

\bibitem{Zeidler1990}
Zeidler, E.: Nonlinear functional analysis and its applications. {II}/{B}.
\newblock Springer-Verlag, New York (1990)

\end{thebibliography}
\end{document}